\documentstyle[12pt,amssymb]{amsart}

\title[Quantum Symmetry Groups of Finite Spaces]
{\Large \bf Quantum Symmetry Groups of Finite Spaces }
\author[Shuzhou Wang]
{\bf Shuzhou Wang}
\dedicatory{Dedicated to Marc A. Rieffel \\
on the occasion of his sixtieth birthday}
\address{Department of Mathematics, University of California,
Berkeley, CA 94720
\newline \indent
Fax: 510-642-8204
}
\email{szwang@@math.berkeley.edu}

\date{}

\newtheorem{DF}{Definition}[section]
\newtheorem{LM}[DF]{Lemma}
\newtheorem{PROP}[DF]{Proposition}
\newtheorem{TH}[DF]{Theorem}
\newtheorem{COR}[DF]{Corollary}
\newtheorem{RMK}[DF]{Remark}
\newtheorem{RMKS}[DF]{Remarks}
\newtheorem{PROB}[DF]{Problem}

\newcommand{\bgdf}{\begin{DF}}
\newcommand{\nddf}{\end{DF}}
\newcommand{\bglm}{\begin{LM}}
\newcommand{\ndlm}{\end{LM}}
\newcommand{\bgprop}{\begin{PROP}}
\newcommand{\ndprop}{\end{PROP}}
\newcommand{\bgth}{\begin{TH}}
\newcommand{\ndth}{\end{TH}}
\newcommand{\bgcor}{\begin{COR}}
\newcommand{\ndcor}{\end{COR}}
\newcommand{\bgrmk}{\begin{RMK}}
\newcommand{\ndrmk}{\end{RMK}}
\newcommand{\bgrmks}{\begin{RMKS}}
\newcommand{\ndrmks}{\end{RMKS}}
\newcommand{\bgprob}{\begin{PROB}}
\newcommand{\ndprob}{\end{PROB}}
\newcommand{\bgequ}{\begin{equation}}
\newcommand{\ndequ}{\end{equation}}
\newcommand{\bgeq}{\begin{eqnarray}}
\newcommand{\ndeq}{\end{eqnarray}}
\newcommand{\bgeqq}{\begin{eqnarray*}}
\newcommand{\ndeqq}{\end{eqnarray*}}

\newcommand{\QED}{\hfill Q.E.D.} 
\newcommand{\vv}{\vspace{4mm}\\}
\newcommand{\vvv}{\vspace{6mm}\\}
\newcommand{\vvvv}{\vspace{1cm}\\}

\numberwithin{equation}{section}

\newcommand{\dfref}[1]{Definition~\ref{#1}}

\newcommand{\thref}[1]{Theorem~\ref{#1}}
\newcommand{\thmref}[1]{Theorem~\ref{#1}}

\begin{document}

\begin{abstract}
We determine the quantum automorphism groups of finite spaces.
These are compact matrix quantum groups in the sense of Woronowicz.
\end{abstract}

\maketitle

\section{
Introduction}

At Les Houches Summer School on Quantum Symmetries in 1995,
Alain Connes posed the following problem:
{\em What is the quantum automorphism group of a space?}
Here the notion of a space is taken in the sense of
noncommutative geometry \cite{Cn}, hence it can be either commutative or
noncommutative.

To put this problem in a proper context, let us recall that the
notion of a group arises most naturally as
symmetries of various kinds of spaces. As a matter of
fact, this is how the notion of a group was discovered historically.
However, the notion of a quantum group was discovered from
several different points of view
\cite{Kac,KP1,Dr1,Wor2,Wor4,Wor5,Wor6,FRT1},
the most important of which is to view quantum groups as deformations of
ordinary Lie groups or Lie algebras, instead of viewing
them as quantum symmetry objects of noncommutative spaces.
In \cite{Mn1}, an important step was made by Manin in
this latter direction, where quantum groups are described
as quantum symmetry objects of quadratic algebras.
(cf also \cite{Wor2} and the book of Sweedler on Hopf algebras.)

In this paper, we solve the problem above for
finite spaces (viz. finite dimensional $C^*$-algebras).
That is, we explicitly determine the quantum automorphism groups of such
spaces. These spaces do not carry the additional geometric (Riemannian)
structures in the sense of \cite{Cn,Cn4}.
The quantum automorphism groups for the latter geometric finite spaces
can be termed quantum isometry groups. At the end of his book \cite{Cn},
Connes poses the problem of finding a finite quantum symmetry
group for the finite geometric space used in his formulation
of the Standard Model in particle physics.
This problem is clearly related to the problem above he
posed at Les Houches Summer School. We
expect that the results in our paper will be useful for
this problem. As a matter of fact, the quantum symmetry group for
the finite geometric space of \cite{Cn}
should be a quantum subgroup of an appropriate quantum automorphism group
described in this paper. The main difficulty is to find
the {\em natural} quantum finite subgroup of the latter that
deserves to be called the quantum isometry group.

This paper can be viewed as a continuation of the work of Manin
\cite{Mn1} in the sense that the quantum groups we consider here
are also quantum symmetry objects. However, it differs
from the work of Manin in three main aspects.
First, the noncommutative spaces on which Manin
considers symmetries are quadratic algebras
and are infinite; while the spaces on which
we consider symmetries are not quadratic and are finite.
Second, Manin's quantum groups are generated by infinitely
many multiplicative matrices and admits many actions on the
spaces in question, one action for each multiplicative matrix
(for the notion of multiplicative matrices, see Manin \cite{Mn1});
while our quantum groups are generated by a single multiplicative
matrix and they act on the spaces in question in {\em one} natural manner.
Finally, Manin's quantum groups do not give rise to natural structures of
$C^*$-algebras in general (see \cite{W5}); while our quantum groups,
besides having a purely algebraic formulation, are compact
matrix quantum groups in the sense of Woronowicz \cite{Wor5}.
Consequently we need to invoke some basic results of
Woronowicz \cite{Wor5}. Loosely speaking, Manin's quantum
groups are noncompact quantum groups.
But to the best knowledge of the author, it is not known as to how one can
make this precise in the strict sense of Woronowicz \cite{Wor8}.
On the other hand, it is natural to expect that quantum automorphism
groups of finite spaces are compact quantum groups without knowing
their explicit descriptions in this paper.

The ideas in our earlier papers \cite{Wang,W1,W5} on universal quantum
groups play an important role in this paper. Note that finite
spaces are just finite dimensional $C^*$-algebras, no
deformation is involved. Moreover, as in \cite{Wang,W1,W5},
the quantum groups considered in this paper are {\em intrinsic}
objects, not as deformations of groups, so they are different from
the quantum groups obtained by the traditional method of deformations
of Lie groups (cf \cite{Dr1,FRT1,Wor4,Wor6,LS1,Rf8,W4}).

We summarize the contents of this paper. In Sect. 2, we recall some basic
notions concerning actions of quantum groups and define
the notion of a quantum automorphism group of a space. The most
natural way to define a quantum automorphism group is by
categorical method, viz, to define it as a universal
object in a certain category of quantum transformation groups.
Sects. 3, 4, 5 are devoted to explicit determination of
quantum automorphism groups for several categories of
quantum transformation groups of the spaces $X_n$, $M_n({\Bbb C})$,
and $\oplus_{k=1}^m M_{n_k}({\Bbb C})$, respectively.
Though the main idea in the construction of quantum automorphism
groups is the same for each of the spaces $X_n$, $M_n({\Bbb C})$ and
$\oplus_{k=1}^m M_{n_k}({\Bbb C})$, the two special cases
$X_n$ and $M_n$ offers interesting phenomena in their own right.
Hence we deal with them separately
and begin by considering the simplest case $X_n$.
In Sect. 6, using the results of sections 3, 4, 5,
we prove that a finite space has a quantum automorphism group
in the category of {\em all} compact quantum transformation groups if and
only if the finite space is $X_n$, and that
a measured finite space (i.e. a finite space endowed with
a positive functional) always has a quantum automorphism group.
\vv
{\em Convention on terminology:} In the following, we will use
interchangeably both the term compact quantum groups
and the term Woronowicz Hopf $C^*$-algebras. When we say that
$A$ is a compact quantum group, we refer to the underlying
geometric object $G$ of $A=C(G)$; when we say that $A$ is a
Woronowicz Hopf $C^*$-algebra, we refer to the `` function algebra''
algebra (cf \cite{Wang,W1,W4,W5}).
\vv
{\em Notation.}
For every natural number $n$, and every *-algebra $A$,
$M_n(A)$ denotes the *-algebra of $n \times n$ matrix with entries in $A$.
We also use $M_n$ to denote $M_n({\Bbb C})$, where $\Bbb C$
is the algebra of complex numbers.
For every matrix $u = (a_{ij}) \in M_n(A)$, $u^t$ denotes the transpose
of $u$; $\bar{u} = (a_{ij}^*)$ denotes the conjugate matrix of $u$;
$u^* = \bar{u}^t$ denotes the adjoint matrix of
$u$ (this defines the ordinary *-operation on $M_n(A)$).
The symbol $X(A)$ denotes the set of all unital *-homomorphism
from $A$ to $\Bbb C$. Finally,
$X_n = \{ x_1, \cdots, x_n \}$ is the finite space with $n$ letters.

\section{
The notion of quantum automorphism groups}

Part of the problem of Connes mentioned in the introduction is
to make precise the notion of a quantum automorphism group, which
we address in this section. First
recall that the usual automorphism group $Aut(X)$ of a space $X$
consists of the set of all transformations on $X$ that preserve
the structure of $X$. A quantum group is not
a set of transformations in general. Thus a naive imitation of the above
definition of $Aut(X)$ for quantum automorphisms will not
work. However, we recapture the definition of $Aut (X)$ from
the following universal property of $Aut(X)$ in the category of
transformation groups of $X$: If $G$ is any group acting on $X$,
then there is a unique morphism of transformation groups
from $G$ to $Aut(X)$. This motivates our definition
\ref{qaut} of quantum automorphism groups below.

The automorphism groups of finite spaces are compact Lie groups
(e.g. $Aut(X_n) = S_n$, the symmetric group on $n$ letters,
and $Aut(M_n) = SU(n)$). For this reason, it is natural to expect
that the quantum automorphism groups of such spaces are compact
quantum groups, viz., Woronowicz Hopf $C^*$-algebras.
We will consider only such quantum groups in this paper.
For basic notions on compact quantum
groups, we refer the reader to \cite{Wor5,Wang,W1}.
Note that for every compact quantum group, there corresponds a
full Woronowicz Hopf $C^*$-algebra and a reduced Woronowicz Hopf
$C^*$-algebra \cite{BS2,W3}. We will assume that all
the Woronowicz Hopf $C^*$-algebras in this paper are full,
as morphisms behave well only with such algebras
(see the discussions in III.7 of \cite{W3}).
Let $A$ be a compact quantum group. Let $\epsilon$ be the unit
of this quantum group (or counit of the full Woronowicz Hopf $C^*$-algebra).
Let $\cal A$ denote the canonical dense Hopf *-subalgebra of $A$
consisting of coefficients of finite dimensional representations of the
quantum group $A$.


\bgdf
\label{qact}
(cf \cite{BS2,Boca1,Pod6})
A {\bf left action} of a compact quantum group $A$ on a
$C^*$-algebra $B$ is a unital *-homomorphism $\alpha$ from $B$ to
$B \otimes A$ such that

(1). $( id_B \otimes \Phi ) \alpha = ( \alpha \otimes id_A ) \alpha,$
where $\Phi$ is the coproduct on $A$;

(2). $(id_B \otimes \epsilon) \alpha = id_B$;

(3). There is a dense *-subalgebra $\cal B$ of $B$, such that $\alpha$
restricts to a right coaction of the Hopf *-algebra $\cal A$ on $\cal B$.

We also call $(A, \alpha)$ a {\bf left quantum transformation group} of $B$.
Let $(\tilde{A}, \tilde{\alpha})$ be another left quantum transformation
group of $B$. We define a {\bf morphism} from $(\tilde{A}, \tilde{\alpha})$
to $(A, \alpha)$ to be a morphism $\pi$ of quantum groups from $\tilde{A}$
to $A$ (which is the same thing as a morphism of Woronowicz Hopf
$C^*$-algebras from $A$ to $\tilde{A}$, see \cite{W1}), such that
$$\tilde{\alpha} = (id_B \otimes \pi) \alpha.$$
It is easy to see that left quantum transformation groups of $B$
form a category with the morphisms defined above.
We call it the {\bf category of left quantum transformation groups} of $B$.
\nddf

Our definition of an action of a quantum group above appears to be different
from the one in \cite{Pod6}, but it is equivalent to the latter.
More precisely, conditions (2) and (3) above are equivalent to
the following density requirement, which is used in
\cite{BS2,Boca1,Pod6} for the definition of an action:
$$(I \otimes A) \alpha (B) \; \; \text{is norm dense in} \; \; B \otimes A,$$
but they are more natural and convenient for our purposes.
It is not clear whether the injectivity condition on $\alpha$
imposed in \cite{BS2,Boca1} is implied by the three conditions
in the definition above.
Our definition coincides with the
notion of actions of groups on spaces when the quantum group $A$ is a
group and $B$ is an ordinary space (simply by reversing the arrows).

The above definition is commonly called the {\em right coaction} of a
unital Hopf $C^*$-algebra.
Note that for the Hopf $C^*$-algebra $A=C(G)$ of continuous
functions over a compact group $G$, the
notion of right coaction of $A$ corresponds to
the notion of left action of $G$ on a $C^*$-algebra $B$.
For this reason, when we are dealing with a compact quantum group $A$,
we call a right coaction of the underlying Woronowicz Hopf
$C^*$-algebra of $A$ a {\bf left action} of the quantum group $A$.
In the following, we will omit the word {\bf left} for actions of
quantum transformation groups. This should not cause confusion.

\bgdf
\label{qfix}
Let $(A, \alpha)$ be a quantum transformation group of $B$.
An element $b$ of $B$ is
said to be {\bf fixed under $\alpha$} (or {\bf invariant under $\alpha$}) if
$$\alpha (b) = b \otimes 1_A.$$
The {\bf fixed point algebra} $A^\alpha$ of the action $\alpha$ is
$$\{ b \in B \; | \; \alpha (b) = b \otimes 1_A \}. $$
The quantum transformation group $(A, \alpha)$ is said to be {\bf ergodic}
if $A^\alpha = {\Bbb C} I$.
A (continuous) functional $\phi$ on $B$ is said to be
{\bf invariant under $\alpha$} if
$$ (\phi \otimes id_A) \alpha (b) = \phi (b) I_A$$
for all $b \in B$.
For a given functional $\phi$ on $B$, we define the {\bf category of quantum
transformation groups of the pair} $(B, \phi)$ to be the category with
objects that leave invariant the functional $\phi$. This is a subcategory of
the category of all quantum transformation groups.
\nddf
Besides the two categories of quantum transformation groups mentioned above,
we also have the category of quantum transformation groups of Kac type for
$B$, which is a full subcategory of the category of
quantum transformation groups of $B$.

\bgdf
\label{qaut}
Let ${\cal C}$ be a category of quantum transformation groups of $B$.
The {\bf quantum automorphism group} of $B$ in ${\cal C}$ is
%
%
a universal final object in the category ${\cal C}$. That is,
if $(\tilde{A}, \tilde{\alpha})$ is an object in this category,
then there is a unique morphism $\pi$ of quantum transformation groups
from $(\tilde{A}, \tilde{\alpha})$ to $(A, \alpha)$.

Let $\phi$ be a continuous functional on the algebra $B$.
We define {\bf quantum automorphism group of the pair} $(B, \phi)$ to be
the universal object in the category of quantum transformation groups
of the pair $(B, \phi)$ (cf \dfref{qact}).
\nddf

From categorical abstract nonsense, the quantum automorphism group of $B$
(in a given category) is unique (up to isomorphism) if it exists.
We emphasize in particular that the notion of a quantum automorphism group
depends on the category of quantum transformation groups of $B$,
not only on $B$. As a matter of fact, for a finite space $B$ other than $X_n$,
we will show in \thmref{mainthm} that the quantum automorphism
group does not exist for the category of all quantum transformation
groups. In the subcategory of quantum
transformation groups of $B$ with objects consisting of compact
transformation groups, the universal object is precisely
the ordinary automorphism group $Aut(B)$, as mentioned in the
beginning of this section.

We will also use the following notion, which generalizes
the usual notion of a faithful group action.

\bgdf
\label{faithful}
Let $(A, \alpha)$ be a quantum transformation group of $B$.
We say that the action $\alpha$ is {\bf faithful} if there is no proper
Woronowicz Hopf $C^*$-subalgebra $A_1$ of $A$ such that
$\alpha$ is an action of $A_1$ on $B$.
\nddf

If $(A, \alpha)$ is a quantum automorphism group in some category of
quantum transformation groups on $B$, then the action $\alpha$
is faithful. We leave the verification of this to the reader as
an exercise.


\section{
Quantum automorphism group of finite space $X_n$}

By the Gelfand-Naimark theorem, we can identify
$X_n = \{ x_1, \cdots, x_n \}$ with the $C^*$-algebra $B=C(X_n)$ of
continuous functions on $X_n$.
The algebra $B$ has the following presentation,
$$
B = C^* \{ e_i \; | \;
e_i^2 = e_i = e_i^*, \; \;
\sum_{r=1}^{n} e_r = 1, \; \; i = 1, \cdots, n \}.
$$
The ordinary automorphism group $Aut(X_n) = Aut(B)$ of
$X_n$ is the symmetric group $S_n$ on $n$ symbols.
We can put the group $S_n$ in the framework of Woronowicz as follows.
As a transformation group, $S_n$ can be thought of as the collection
of all permutation matrices
$$ g =
\left(
\begin{array}{cccc}
a_{11} & a_{12} & \cdots & a_{1n}    \\
a_{21} & a_{22} & \cdots & a_{2n}    \\
\cdots  & \cdots  & \cdots & \cdots     \\
a_{n1} & a_{n2} & \cdots & a_{nn}
\end{array}
\right).
 $$
 When $g$ varies in $S_n$, the $a_{ij}$'s ($i,j = 1, \cdots, n$)
 are functions on the group $S_n$ satisfying the following relations:
\bgeq
\label{r1}  a_{ij}^2 = a_{ij} = a_{ij}^*,
\; \; \; i,j = 1, \cdots, n, \\
\label{r2}  \sum_{j = 1}^{n} a_{ij} = 1 ,
\; \; \; i = 1, \cdots, n, \\
\label{r3}  \sum_{i = 1}^{n} a_{ij} = 1 ,
\; \; \; i = 1, \cdots, n.
\ndeq
It is easy to see that the commutative $C^*$-algebra generated by
the above commutation relations is the Woronowicz Hopf $C^*$-algebra
$C(S_n)$. In other words, the group $S_n$ is completely determined by
these relations. The following theorem shows that we have obtained much more:
If we remove the condition that the $a_{ij}$'s commute with each other,
these relations defines the quantum automorphism group of $X_n$.

\bgth
\label{thm1}
Let $A$ be the $C^*$-algebra with generators $a_{ij}$ ($i,j = 1, \cdots, n$)
and defining relations \eqref{r1} -- \eqref{r3}. Then

(1). $A$ is a compact quantum group of Kac type;

(2). The formulas
$$
\alpha (e_j) = \sum_{i=1}^{n} e_i \otimes {a}_{ij}, \; \; \; j = 1, \cdots, n
$$
defines a quantum transformation group $(A, \alpha)$ of $B$. 
It is the quantum automorphism group of $B$ 
in the category of all compact quantum transformation groups
(hence also in the category of compact quantum groups of Kac type)
of $B$, 
and it contains the ordinary automorphism group $Aut(X_n) = S_n$ (in fact,
$\{(\chi (a_{ij})) \; | \; \chi \in X(A) \}$ is precisely the set of
permutation matrices).
\ndth

Because of (2) above, we will denote the quantum group above by
$A_{aut} (X_n)$. We will call it the {\bf quantum permutation group}
on $n$ symbols.

\pf (1). 
It is easy to check that there is a well-defined homomorphism $\Phi$ from $A$
to $A \otimes A$ with the property
$$\Phi (a_{ij}) = \sum_{k=1}^n a_{ik} \otimes a_{kj},
\; \; \; i,j = 1, \cdots, n.$$
Using \eqref{r1} -- \eqref{r3},
it is also easy to check that $u=(a_{ij})$ is an orthogonal matrix.
Hence $(A, u)$ is a quantum subgroup of $A_o(n)$, so it is of Kac type
(cf \cite{Wang,W1,W5}).

To prove (2),
note that the generators $\{ e_i \}_{i=1}^n$ form a basis of the vector
space $B$, so an action $\tilde{\alpha}$ of any quantum group
$\tilde{A}$ on $B$ is uniquely determined by its effect on the $e_i$'s:
$$
\tilde{\alpha} (e_j) = \sum_{i=1}^{n} e_i \otimes \tilde{a}_{ij},
\; \; \; j = 1, \cdots, n.
$$
%
%

The condition that $\tilde{\alpha}$ is a *-homomorphism together
with the equations
$$e_i^2 = e_i = e_i^* ,
\; \; \; i = 1, \cdots, n \\
$$
shows that the $\tilde{a}_{ij}$'s satisfy the relations \eqref{r1}.
The condition that $\tilde{\alpha}$ is a unital homomorphism together
with the equation
$$ \sum_{i=1}^{n} e_i = 1 $$
shows that the $\tilde{a}_{ij}$'s satisfy \eqref{r2}.
Let $\tilde{u} = (\tilde{a}_{ij})$. Then we have
$$\tilde{u} \tilde{u}^* = I_n.$$
The condition in \dfref{qact}.(2) means that
$$\epsilon (\tilde{a}_{ij}) = \delta_{ij}, \; \; \; i,j = 1, \cdots, n.$$
By condition (3) of \dfref{qact}, the $\tilde{a}_{ij}$'s are
in $\tilde{\cal A}$. Hence by Proposition 3.2 of \cite{Wor5},
it follows that $\tilde{u} = (\tilde{a}_{ij})$ is a non-degenerate smooth
representation of the quantum group $\tilde{A}$.
In particular, $\tilde{u}$ is also left invertible,
$$\tilde{u}^* \tilde{u} = I_n.$$
This implies that the $\tilde{a}_{ij}$'s satisfy the
relations \eqref{r3}.
From these we see that $(A, \alpha)$ is a universal
quantum transformation group of
$B$: there is a unique morphism $\pi$ of quantum
transformation groups from $(\tilde{A}, \tilde{\alpha})$ to
$(A, \alpha)$ such that
$$\pi (a_{ij}) = \tilde{a}_{ij} ,
\; \; \; i,j = 1, \cdots, n.$$

It is clear that the maximal
subgroup of the quantum group $A$ is $S_n$, that is, the set
$\{(\chi (a_{ij})) \; | \; \chi \in X(A) \}$ is precisely the set of
permutation matrices.
\QED
\vv
{\em Remarks.}
(1).
For each pair $i,j$, let $A_{ij}$ be the group $C^*$-algebra
$ C^* ( {\Bbb Z} / 2{\Bbb Z} ) $ with generator $p_{ij}$,
$p_{ij}^2 = p_{ij} = p_{ij}^*$ ($i,j = 1, \cdots, n$).
Then the $C^*$-algebra $A$ is isomorphic to the following quotient
$C^*$-algebra of the free product of the $A_{ij}$'s:
$$( \ast_{i,j=1}^n A_{ij} ) /
< \sum_{r = 1}^{n} p_{rj} = 1 =
\sum_{s = 1}^{n} p_{is}, \; \; i,j = 1, \cdots, n >.$$
From this we see that for $n \leq 3$, $A = C(S_n)$, for $n \geq 4$, $A$
is noncommutative and infinite dimensional.

(2).
Let $\phi$ be the unique $S_n$-invariant probability measure on $X_n$.
Then it is easy to see that $\phi$ is a fixed functional under the action of
the quantum group $A_{aut}(X_n)$ defined in \thref{thm1}. Hence
$A_{aut}(X_n)$ is also the quantum automorphism group for the
pair $(X_n, \phi)$.

(3).
Let $Q > 0$ be a positive $n \times n$ matrix.
Let $A^Q_{aut}(X_n)$ be the $C^*$-algebra with generators
$a_{ij}$ ($i,j = 1, \cdots, n$) and the defining relations given by
\eqref{r1} -- \eqref{r2} along with the following set of relations:
\bgeq
\label{r4}  u^t Q u Q^{-1} = I_n =  Q u Q^{-1} u^t ,
\ndeq
where $u = (a_{ij})$.
Then it not hard to verify that $(A^Q_{aut}(X_n), \alpha)$ is
a compact quantum transformation subgroup of the one defined in
\thref{thm1} (hence the $a_{ij}$'s also satisfy the relations \eqref{r3}),
here $\alpha$ is as in \thref{thm1}.
Note also for $Q = I_n$, $A^Q_{aut}(X_n) = A_{aut}(X_n)$.
%
%

\section{
Quantum automorphism group of finite space $M_n ({\Bbb C})$}

\noindent
{\em Notation.} Let $u = (a^{kl}_{ij})_{i,j,k,l=1}^n$ and
$v = (b^{kl}_{ij})_{i,j,k,l=1}^n$ with entries from a *-algebra.
Define $ u v$ to be the
matrix whose entries are given by
$$(u v)^{kl}_{ij} = \sum_{r,s=1}^n a^{kl}_{rs} b^{rs}_{ij},
\; \; \; i,j,k,l = 1, \cdots, n.$$

Let $\psi = Tr$ be the trace functional on $M_n$ (so
$\phi = \frac{1}{n} \psi$ is the unique $Aut(M_n)$-invariant state on
$M_n$). The $C^*$-algebra $M_n$ has the following presentation
$$
B = C^* \{ e_{ij} \; | \;
e_{ij} e_{kl} = \delta_{jk} e_{il} , \;
e_{ij}^* = e_{ji},  \; \;
\sum_{r=1}^n e_{rr} = 1, \; \; i,j,k,l = 1, \cdots, n \} .
$$

\bgth
\label{thm3}
Let $A$ be the $C^*$-algebra with generators $a^{kl}_{ij}$
and the following defining
relations \eqref{r5} -- \eqref{r9}:
\bgeq
\label{r5}  \sum_{v=1}^n a^{kv}_{ij} a^{vl}_{rs} = \delta_{jr} a^{kl}_{is},
\; \; \; i,j,k,l,r,s = 1, \cdots, n, \\
\label{r6}  \sum_{v=1}^n a^{sr}_{lv} a^{ji}_{vk} = \delta_{jr} a^{si}_{lk},
\; \; \; i,j,k,l,r,s = 1, \cdots, n,   \\
\label{r7}  {a^{kl}_{ij}}^* = a^{lk}_{ji},
\; \; \; i,j,k,l = 1, \cdots, n, \\
\label{r8}  \sum_{r=1}^n a^{kl}_{rr} = \delta_{kl},
\; \; \; k,l = 1, \cdots, n, \\
\label{r9}  \sum_{r=1}^n a^{rr}_{kl} = \delta_{kl},
\; \; \; k,l = 1, \cdots, n.
\ndeq
Then

(1). $A$ is a compact quantum group of Kac type;

(2). The formulas
$$
\alpha (e_{ij}) = \sum_{k,l=1}^{n} e_{kl} \otimes {a}^{kl}_{ij},
\; \; \; \; i,j = 1, \cdots, n
$$
defines a quantum transformation group $(A, \alpha)$ of $(M_n, \psi)$.
It is the quantum automorphism
group of $(M_n, \psi)$ in the category of compact quantum transformation
groups (hence also in the category of compact quantum groups of Kac type)
of $(M_n, \psi)$,
and it contains the ordinary automorphism group $Aut(M_n) = SU(n)$.
\ndth

We will denote the quantum group above by
$A_{aut} (M_n)$.

\pf
(1).
It is easy to check that the matrix $u = (a^{kl}_{ij})$ as well
as its conjugate $\bar{u} = ({a^{kl}_{ij}}^*)$ are both unitary
matrices, and that the formulas
$$\Phi ( a^{kl}_{ij} ) = \sum_{r,s=1}^n a^{kl}_{rs} \otimes a^{rs}_{ij},
\; \; \; i,j,k,l = 1, \cdots, n $$
gives a well-defined map from $A$ to $A \otimes A$ (this is the coproduct).
Hence $A$ is a quantum subgroup of $A_u(m)$ (with $m = n^2$),
so it is of Kac type (cf \cite{Wang,W1,W5}).

(2).
Let $(\tilde{A}, \tilde{\alpha})$ be any quantum transformation group of
$M_n$. Being a basis for the vector space $M_n$, the $e_{ij}$'s uniquely
determine the action $\tilde{\alpha}$:
$$
\tilde{\alpha} (e_{ij}) = \sum_{k,l=1}^{n} e_{kl} \otimes \tilde{a}^{kl}_{ij},
\; \; \; i,j = 1, \cdots, n.
$$

The condition that $\tilde{\alpha}$ is a homomorphism together
with the equations
$$
e_{ij} e_{kl} = \delta_{jk} e_{il} , \;  \; \; i,j,k,l = 1, \cdots, n
$$
shows that the $\tilde{a}^{kl}_{ij}$'s satisfy \eqref{r5}.
The condition that $\tilde{\alpha}$  preserves the *-operation together
with the equations
$$
e_{ij}^* = e_{ji}, \; \; \; i,j = 1, \cdots, n
$$
shows that the $\tilde{a}^{kl}_{ij}$'s satisfy \eqref{r7}.
The condition that $\tilde{\alpha}$ preserves the units together
with the identity
$$
\sum_r e_{rr} = 1
$$
shows that the $\tilde{a}^{kl}_{ij}$'s satisfy \eqref{r8}.
The condition that $\tilde{\alpha}$ leaves the trace $\psi$
invariant shows that the $\tilde{a}^{kl}_{ij}$'s satisfy \eqref{r9}.

To show that the $\tilde{a}^{kl}_{ij}$'s satisfy \eqref{r6}, first
it is an easy check that
$$ \tilde{u}^* \tilde{u} = I_n^{\otimes 2},$$
where $\tilde{u} = (\tilde{a}^{kl}_{ij})_{i,j,k,l=1}^n$.
By condition (3) of \dfref{qact}, the $\tilde{a}^{kl}_{ij}$'s are
in $\tilde{\cal A}$. Hence by Proposition 3.2 of \cite{Wor5}, we see that
$\tilde{u}$ is a non-degenerate smooth
representation of the quantum group $\tilde{A}$. In
particular, $\tilde{u}$ is also right invertible,
$$ \tilde{u} \tilde{u}^* = I_n^{\otimes 2}, $$
which means that
$$\sum_{i,j=1}^n \tilde{a}^{kl}_{ij} \tilde{a}^{sr}_{ji} = \delta_{kr}
\delta_{ls}, \; \; \; k,l,r,s = 1, \cdots, n.$$
From these relations and the relations \eqref{r5}, \eqref{r7}-\eqref{r9}, we
deduce that both matrices $\tilde{u}$ and $\tilde{u}^t$
are unitary. This shows that the quantum group $A_1$
generated by the coefficients $\tilde{a}^{kl}_{ij}$
is a compact quantum group of Kac type.
That is, the antipode $\tilde{\kappa}$ is a bounded
*-antihomomorphism when restricted to $A_1$. Put
$$v = (b^{kl}_{ij}) = (\tilde{\kappa} (\tilde{a}^{kl}_{ij}))
= (\tilde{a}^{ji}_{lk}).$$
Then in the opposite algebra ${A_1}^{op}$ (which has the
same elements as $A_1$ with multiplication reserved),
the $b^{kl}_{ij}$'s satisfy the relations \eqref{r5}, which means that
the $\tilde{a}^{kl}_{ij}$'s  satisfy
the relations \eqref{r6} in the algebra $\tilde A$.

From the above consideration we see that $(A, \alpha)$ is a quantum
transformation group of $M_n$, and that there is a unique morphism $\pi$ of
quantum groups from $\tilde{A}$ to $A$ such that
$$\pi (a^{kl}_{ij}) = \tilde{a}^{kl}_{ij} ,
\; \; \; i,j,k,l = 1, \cdots, n.$$
It is routine to check that $\pi$ is the unique morphism $\pi$ of quantum
transformation groups from $(\tilde{A}, \tilde{\alpha})$ to
$(A, \alpha)$.

From the relations \eqref{r5} -- \eqref{r9}, one can show that each matrix
$(\chi (a^{kl}_{ij}))$ ($\chi \in X(A_{aut}(M_n))$) defines
an automorphism of $M_n$ by the formulas in \thref{thm3}.(2).
This means that the maximal subgroup $X(A_{aut}(M_n))$
is naturally embedded in $Aut(M_n)$.
Conversely, it is clear that $Aut(M_n)$ can be embedded as a subgroup of the
maximal subgroup $X(A_{aut}(M_n))$ of $A_{aut}(M_n)$.
\QED
\vv
{\em Remark.}
Consider the quantum group $(A_u(n), (a_{ij}))$ (cf \cite{W1,W5}).
Put $\tilde{a}^{kl}_{ij} = a_{ki} a_{lj}^*$. Then the
$\tilde{a}^{kl}_{ij}$'s satisfies the relations \eqref{r5} -- \eqref{r9}.
From this we see that the $\tilde{a}^{kl}_{ij}$'s determines a quantum
subgroup of $A_{aut} (M_n)$. Hence the Woronowicz Hopf $C^*$-algebra
$A_{aut}(M_n)$ is noncommutative and nococommutative.
How big is the subalgebra of $A_u(n)$ generated by the $\tilde{a}^{kl}_{ij}$?
An answer to this question will shed light on the structure of
the $C^*$-algebra $A_{aut}(M_n)$.

\bgprop
\label{thm4}
Let $Q > 0 $ be a positive matrix in $M_n ({\Bbb C}) \otimes M_n ({\Bbb C})$.
Let $A$ be the $C^*$-algebra with generators $a^{kl}_{ij}$ and defining
relations given by \eqref{r5}, \eqref{r7}, \eqref{r8}, along with the
following set of relations:
\bgeq
\label{r10} u^* Q u Q^{-1} = I_n^{\otimes 2} =  Q u Q^{-1} u^* ,
\ndeq
where $u = (a^{kl}_{ij})$. Then $A$ is 
a compact quantum group that acts faithfully on $M_n$ in the following manner,
$$
\alpha (e_{ij}) = \sum_{k,l=1}^n e_{kl} \otimes a^{kl}_{ij},
\; \; \; i,j = 1, \cdots, n, $$
and its maximal subgroup is isomorphic to a subgroup
of $Aut(M_n) \cong SU(n)$.
Any faithful compact quantum transformation group of $M_n$ is a quantum
subgroup of $(A, \alpha)$ for some positive $Q$.
\ndprop
\pf
First we show that $A$ is a compact quantum group.
Let $v = Q^{1/2} u Q^{- 1/2}$. Then \eqref{r10} is equivalent to
$$ v^* v = I_n^{\otimes 2} = v v^*.$$
Hence the $C^*$-algebra $A$ is well defined.
The set of relations in \eqref{r10} shows that $u$ is invertible.
We claim that $u^t$ is also invertible. For simplicity of
notation in the following computation, let $\tilde{Q} =
(\tilde{q}^{kl}_{ij}) = Q^{-1}$. Then \eqref{r10} becomes
$$
\sum_{k,l,r,s,x,y = 1}^n
a^{lk}_{ij} q^{kl}_{rs} a^{rs}_{xy} \tilde{q}^{xy}_{ef}
= \delta_{ef}^{ij} =
\sum_{k,l,r,s,x,y = 1}^n
q^{ij}_{kl} a^{kl}_{rs} \tilde{q}^{rs}_{xy} a^{yx}_{fe},
$$
where $i,j,e,f = 1, \cdots, n.$
Put $P = (p^{kl}_{ij})$ and $\tilde{P} = (\tilde{p}^{kl}_{ij})$, where
$$p^{kl}_{ij} = q^{lk}_{ij}, \; \; \;
\tilde{p}^{kl}_{ij} = q^{kl}_{ji}, \; \; \; i,j,k,l = 1, \cdots, n.
$$
Then $P^{-1} = \tilde{P}$, and the relations \eqref{r10}
becomes
$$
 u^t P u P^{-1} = I_n^{\otimes 2} =  P u P^{-1} u^t .
$$
This proves our claim.

Now it is easy to check that $A$ is a compact matrix quantum group
with coproduct $\Phi$ given by the same formulas as in the
proof of \thmref{thm3}.(1).

Let $(\tilde{A}, \tilde{\alpha})$ be a faithful quantum transformation group
of $M_n$. We saw in the proof of \thmref{thm3} that there are
elements $\tilde{a}^{kl}_{ij}$ ($i,j,k,l = 1, \cdots, n$) in the $C^*$-algebra
$\tilde{A}$ that satisfy the relations \eqref{r5}, \eqref{r7} and \eqref{r8}.
The condition in \dfref{qact}.(2) means that
$$\epsilon (\tilde{a}^{kl}_{ij}) = \delta^{kl}_{ij},
\; \; \; i,j,k,l = 1, \cdots, n.$$
By condition (3) of \dfref{qact}, the $\tilde{a}_{ij}$'s are
in $\tilde{\cal A}$.
Hence by Proposition 3.2 of \cite{Wor5}, this implies that
$\tilde{u} = (\tilde{a}^{kl}_{ij})$ is a non-degenerate smooth
representation of the quantum group $\cal A$.
From the proof of Theorem 5.2 of \cite{Wor5}, with
$$Q = (id \otimes \tilde{h}) (\tilde{u}^* \tilde{u}),$$
we have $Q > 0$ and $\tilde{u}$ satisfies \eqref{r10}.
The assumption that $(\tilde{A}, \tilde{\alpha})$ is faithful implies
that $\tilde A$ is generated by the elements
$\tilde{a}^{kl}_{ij}$ ($i,j = 1, \cdots, n$). This shows that $(A, \alpha)$
is a well defined faithful quantum transformation group of $M_n$
and that the compact quantum transformation group
$({\tilde A}, {\tilde \alpha})$ is a quantum subgroup of $(A, \alpha)$.

Let $\chi \in X(A)$. From the defining relations for $A$, we see
that $(\chi(a_{kl,ij}))$ defines an ordinary transformation
for $M_n$ via the formulas in \thmref{thm4}. Hence
the maximal subgroup $X(A)$ is embedded in $Aut(M_n)$.
\QED
\vv
{\em Note.}
We will denote the quantum group above by $A_{aut}^Q (M_n)$.
If $Q = I_n^{\otimes 2}$, then it is easy to see that
the square of the coinverse (i.e. antipode) map is the identity map.
From this one can show that this quantum group reduces to the quantum
group $A_{aut}(M_n)$ in \thmref{thm3}.

\section{
Quantum automorphism group of finite space
$\bigoplus_{k=1}^m M_{n_k} ({\Bbb C})$}

\noindent
{\em Notation.} Let $u = (a^{kl}_{rs, xy})$ and
$v = (b^{kl}_{rs, xy})$ be two matrices with entries from a *-algebra, where
$$
\; \; k,l = 1, \cdots, n_x,
\; \; r,s = 1, \cdots, n_y,
\; \; x,y = 1, \cdots, m.
$$
Define $ u v$ to be the
matrix whose entries are given by
$$(u v)^{kl}_{rs, xy} = \sum_{p=1}^m
\sum_{i,j=1}^{n_p} a^{kl}_{ij, xp} b^{ij}_{rs, py} .$$

Using the same method as above, we now study the quantum
automorphism group of the finite space $B = \bigoplus_{k=1}^m M_{n_k}$, where
$n_k$ is a positive integer. The $C^*$-algebra $B$ has the following
presentation
\bgeqq
B = C^* \{ e_{kl,i} \; |
& & e_{kl,i} e_{rs,j} = \delta_{ij} \delta_{lr} e_{ks} , \;
    e_{kl,i}^* = e_{lk,i},  \; \;
    \sum_{q=1}^m \sum_{p=1}^{n_q} e_{pp,q} = 1, \\
& & k, l = 1, \cdots, n_i , \; \;
    r, s = 1, \cdots, n_j , \; \;
    i, j = 1, \cdots, m
    \} .
\ndeqq
Let $\psi$ be the positive functional on $B$ defined by
$$\psi (e_{kl,i}) = Tr(e_{kl, i}) = \delta_{kl}, \; \;
k, l = 1, \cdots, n_i, \; \;
i = 1, \cdots, m.
$$
The defining relations for the quantum group of $(B, \psi)$
are obtained as a combination of the relations
of the quantum automorphism groups $A_{aut}(X_n)$ and $A_{aut}(M_n)$.

\bgth
\label{thm5}
Let $A$ be the $C^*$-algebra with generators $a^{kl}_{rs, xy}$
$$
\; \; k,l = 1, \cdots, n_x,
\; \; r,s = 1, \cdots, n_y,
\; \; x,y = 1, \cdots, m,
$$
and the following defining relations \eqref{r11} -- \eqref{r15}:
\bgeq
\label{r11}  \sum_{v=1}^{n_x} a^{kv}_{ij, xy} a^{vl}_{rs, xz}
= \delta_{jr} \delta_{yz} a^{kl}_{is, xy},
\ndeq
$$
\; \; i,j = 1, \cdots, n_y,
\; \; r,s = 1, \cdots, n_z,
\; \; k,l = 1, \cdots, n_x,
\; \; x,y,z = 1, \cdots, m,
$$
\bgeq
\label{r12}  \sum_{v=1}^{n_x} a^{sr}_{lv, yx} a^{ji}_{vk, zx}
= \delta_{jr} \delta_{yz} a^{si}_{lk, yx},
\ndeq
$$
\; \; i,j = 1, \cdots, n_z,
\; \; r,s = 1, \cdots, n_y,
\; \; k,l = 1, \cdots, n_x,
\; \; x,y,z = 1, \cdots, m,
$$
\bgeq
\label{r13}  {a^{kl}_{ij, yz}}^* = a^{lk}_{ji, yz},
\ndeq
$$
\; \; i,j = 1, \cdots, n_z,
\; \; k,l = 1, \cdots, n_y,
\; \; y,z = 1, \cdots, m,
$$
\bgeq
\label{r14}  \sum_{z=1}^m \sum_{r=1}^{n_z} a^{kl}_{rr, yz} = \delta_{kl},
\; \; \; k,l = 1, \cdots, n_y,
\; \; \; y = 1, \cdots, m, \\
\label{r15}  \sum_{y=1}^m \sum_{r=1}^{n_y} a^{rr}_{kl, yz} = \delta_{kl},
\; \; \; k,l = 1, \cdots, n_z,
\; \; \; z = 1, \cdots, m.
\ndeq
Then

(1). $A$ is a compact quantum group of Kac type;

(2). The formulas
$$
\alpha (e_{rs, j}) =
\sum_{i=1}^m \sum_{k,l}^{n_i} e_{kl, i} \otimes {a}^{kl}_{rs, ij} ,
\; \; \; r,s = 1, \cdots, n_j, \; \; \; j = 1, \cdots, m
$$
defines a quantum transformation group $(A, \alpha)$ of $(B, \psi)$.
This is the quantum automorphism group of $(B, \psi)$
in the category of compact quantum transformation
groups (hence also in the category of compact quantum groups of Kac type)
of $(B, \psi)$, and it contains the ordinary automorphism group $Aut(B)$.
\ndth

We will denote the quantum group above by $A_{aut} (B)$.

\pf
The proof of this theorem follows the lines of the proof of
\thref{thm3}. The coproduct is given by
$$\Phi ( a^{kl}_{ij, xy} ) =
\sum_{p=1}^m \sum_{r,s=1}^{n_p} a^{kl}_{rs, xp} \otimes a^{rs}_{ij, py},
\; \; \; k, l = 1, \cdots, n_x,
\; \; \; x, y = 1, \cdots, m.
$$
\QED

Note that when $n_k = 1$ for all $k$, then the quantum group $A_{aut}(B)$
reduces to the quantum group $A_{aut}(X_n)$ in \thref{thm1}, and
when $m = 1$, $A_{aut}(B)$ reduces to the quantum group $A_{aut}(M_n)$
in \thref{thm3}.

Let $Q = (q^{kl}_{rs,xy}) > 0 $
($ k, l = 1, \cdots, n_x, \; \;
r, s = 1, \cdots, n_y, \; \;
x, y = 1, \cdots, m $)
be a positive matrix with complex entries.
Define $\delta^{kl}_{rs,xy} $ to be $1$ if $k=r, l=s, x=y$ and
$0$ otherwise, and let $I$ be the matrix with entries
 $\delta^{kl}_{rs,xy} $,
where
$$ k, l = 1, \cdots, n_x, \; \;
   r, s = 1, \cdots, n_y, \; \;
   x, y = 1, \cdots, m .$$

\bgprop
\label{thm6}
Let $Q$ and $I$ be as above.
Let $A$ be the $C^*$-algebra with generators $a^{kl}_{rs, xy}$
$$
\; \; k,l = 1, \cdots, n_x,
\; \; r,s = 1, \cdots, n_y,
\; \; x,y = 1, \cdots, m,
$$
and defining relations \eqref{r11}, \eqref{r13}, \eqref{r14},
along with the following set of relations:
\bgeq
\label{r16} u^* Q u Q^{-1} = I =  Q u Q^{-1} u^* ,
\ndeq
where $u = (a^{kl}_{rs, xy})$. Then $A$ is
a compact quantum group that acts faithfully on $B$ in the following manner,
$$
\alpha (e_{rs, j}) =
\sum_{i=1}^m \sum_{k,l}^{n_i} e_{kl, i} \otimes {a}^{kl}_{rs, ij} ,
\; \; \; r,s = 1, \cdots, n_j, \; \; \; j = 1, \cdots, m.
$$
Any faithful compact quantum transformation group of $B$ is a quantum
subgroup of $(A, \alpha)$ for some positive $Q$.
\ndprop
\pf
The proof follows the lines of \thref{thm4}.
\QED
\vv
We will denote the quantum group above by $A^Q_{aut} (B)$, or simply by
$A^Q_{aut}$. When $Q = I_n^{\otimes 2}$,
then $A^Q_{aut} (B)$ is just $A_{aut} (B)$.
%

Note that for $n_k$'s distinct, the automorphism group
$Aut ( \oplus_{k=1}^m M_{n_k} )$ is isomorphic to the group
$\times_{k=1}^m Aut( M_{n_k} )$. A natural problem related to this is

\bgprob
For $n_k$'s distinct, is the quantum automorphism group
$A_{aut} ( \oplus_{k=1}^m M_{n_k} )$ isomorphic to
the quantum group $\otimes_{k=1}^m A_{aut}( M_{n_k} )$ (cf \cite{W2}).
\ndprob

For each fixed $1 \leq k_0 \leq m$, $A_{aut}(M_{k_0})$ as defined
in the last section is a quantum subgroup of $A_{aut}(B)$.
(This is seen as follows. Let
$\tilde{a}^{kl}_{rs, xy} = \delta_{x k_0} \delta_{y k_0} a^{kl}_{rs}$, where
the $a^{kl}_{rs}$'s are generators of $A_{aut}(M_{n_{k_0}})$. Then
the $\tilde{a}^{kl}_{rs, xy}$'s satisfy the defining relations for
$A_{aut}(B)$.)
Note also that if $n_k = n$ for all $k$, then $A_{aut}(X_m)$ is a quantum
subgroup of $A_{aut}(B)$.
(This is seen as follows. Let
$\tilde{a}^{kl}_{rs, xy} = \delta_{k r} \delta_{l s} a_{xy}$, where
the $a_{xy}$'s are generators of $A_{aut}(X_m)$. Then
the $\tilde{a}^{kl}_{rs, xy}$'s satisfy the defining relations for
$A_{aut}(B)$.)
In view of the fact that the ordinary automorphism
group $Aut(\oplus_1^m M_n)$ is isomorphic to the semi-direct product
$SU(n) \rtimes S_m$, it would be interesting to solve the following problem.

\bgprob
Is it possible to express $A_{aut} (\oplus_1^m M_n)$ in terms
of $A_{aut}(M_n)$ and $A_{aut}(X_m)$ as a certain semi-direct
product that generalizes \cite{W2}?
\ndprob

\section{
The main result}

Summarizing the previous sections, we can now
state the main result of this paper.

\bgth
\label{mainthm}
Let $B$ be a finite space of the form $\oplus_{k=1}^m M_{n_k}$.

(1). Quantum automorphism group of $B$ exists in the
category of (left) quantum transformation groups if and only if $B$ is the
finite space $X_m$.

(2). The quantum automorphism group for $(B, \psi)$ exists and is
defined as in \thref{thm5} (see also \thref{thm1}, \thref{thm3}).
\ndth
\pf
(1).
If $B$ is $X_m$, we saw in \thref{thm1} that $A_{aut}(X_m)$ is the
quantum automorphism group of $X_m$ in the category
of all quantum transformation groups.

Now assume that $B \neq C(X_m)$, and assume that quantum
automorphism group of $B$ exists in the category of all quantum
transformation groups. Call it $(A_0, \alpha_0)$. As in
\thref{thm5} and \thref{thm6}, $\alpha_0$ is determined by
its effect on the basis $e_{rs, j}$ of $B$,
$$
\alpha_0 (e_{rs, j}) =
\sum_{i=1}^m \sum_{k,l}^{n_i} e_{kl, i} \otimes \tilde{a}^{kl}_{rs, ij} ,
\; \; \; r, s = 1, \cdots, n_j, \; \; \; j = 1, \cdots, m.
$$
Since $(A_0, \alpha_0)$ is the quantum automorphism group of $B$, the action
$\alpha_0$ is faithful (cf \dfref{faithful}). This implies that
the $\tilde{a}^{kl}_{rs, ij}$'s generates the $C^*$-algebra
$A_0$. As in \thref{thm6} (see also \thref{thm4}),
there is a positive $Q_0$, such that the $\tilde{a}^{kl}_{rs, xy}$'s
satisfy the relations \eqref{r11}, \eqref{r13}, \eqref{r14},
along with the following set of relations:
\bgeq
\label{r17}
\tilde{u}^* Q_0 \tilde{u} Q_0^{-1} =
I =  Q_0 \tilde{u} Q_0^{-1} \tilde{u}^* ,
\ndeq
where $\tilde{u} = (\tilde{a}^{kl}_{rs, xy})$. By the universal
property of $(A_0, \alpha_0)$, we conclude that $A_0 = A^{Q_0}_{aut}$
(see also the last statement in \thref{thm6}).
For every positive $Q$, the unique morphism from
$(A^Q_{aut}, \alpha)$ to $(A_0, \alpha_0)$ sends the generators
$\tilde{a}^{kl}_{rs, xy}$ of $A^{Q_0}_{aut}$
to the corresponding generators $a^{kl}_{rs, xy}$ of $A^Q_{aut}$
(again because of faithfulness of the quantum transformation group
$A^{Q}_{aut}$ and the universality of $A^{Q_0}_{aut}$).
Hence the generators $a^{kl}_{rs, xy}$ also satisfy
the relations \eqref{r17}. This is impossible because
we can choose $Q$ so that $A^Q_{aut}$ and $A^{Q_0}_{aut}$
have different {\em classical points} in the {\em vector space} with
coordinates $a^{kl}_{rs, xy}$
($ k,l = 1, \cdots, n_x,
\; \; r,s = 1, \cdots, n_y,
\; \; x,y = 1, \cdots, m$).

(2). This is proved in the previous sections.
\QED
\vv
{\em Concluding Remarks.}
(1).
In this paper, we only described the quantum automorphism
group of $(B, \psi)$ for the special choice of functional $\psi$,
because this quantum automorphism group is closest to the
ordinary automorphism group $Aut(B)$ of $B$, and it contains
the latter. One can also use the same method to describe quantum
automorphism groups of $B$ endowed with other functionals
or a collection of functionals.

(2).
For each $1 \leq k \leq n$, consider the delta measure $\chi_k$
on $X_n$ corresponding to the point $x_k$.
Then the quantum automorphism group of $(X_n, \chi_k)$
is isomorphic to the {\em quantum permutation group} of the
space $X_{n-1}$, just as in the case of ordinary permutation groups.

(3).
If we remove condition (3) in \dfref{qact}, then we obtain the notion
of an action of a quantum semi-group on a $C^*$-algebra.
The relations \eqref{r11}, \eqref{r13}, \eqref{r14} define
{\em the universal quantum semi-group} $E(B)$ acting on $B$,
even though $B$ is not a quadratic algebra in the
sense of Manin \cite{Mn1}. From the main theorem of this paper,
the Hopf envelope $H(B)$ of this quantum semi-group in the sense of Manin
cannot be a compact quantum group (see also the last
section of \cite{W5}).

After this paper was submitted for publication, we received the
papers \cite{DHS,DNS}, where a finite quantum group symmetry $A(F)$
for $M_3$ is described, following the work of Connes \cite{Cn4}.
The finite quantum group $A(F)$ in these papers is not
a finite quantum group in the sense of \cite{Wor5} (because
it does not have a compatible $C^*$ norm),
so it cannot be a quantum subgroup of the
COMPACT quantum symmetry groups $A_{aut}(M_3)$ and $A^Q_{aut}(M_3)$
in our paper; but it is a quantum subgroup of the Hopf envelope $H(B)$
of the quantum semi-group mentioned in the last paragraph.

Our paper gives solutions to the ``intricate problem'' mentioned in
the end of section 2 of the paper \cite{DNS}: find the
biggest quantum group acting on $M_3$.
This ``intricate problem'' has two solutions:
the first, \thref{mainthm}, solves the
problem in the category of compact quantum groups; the second, the remarks
in the last two paragraphs, solves the problem in the category of all
quantum groups--Hopf algebras that need not have $C^*$-norms.

(4).
In \cite{Mn1}, the quantum group $SU_q(2)$ is described as
the quantum automorphism group of the quantum plane (i.e. the
deformed plane). In view of the fact that the automorphism group
$Aut(M_2)$ is $SU(2)$, one might be able to describe $SU_q(2)$
as a quantum automorphism group of the non-deformed
space $M_2$ endowed with a collection of functionals.
\vvv
\begin{center}
{ \Large \bf
An Appendix}
\end{center}
\vspace{2mm}
In \cite{W5}, we introduced a compact matrix quantum group
$A_o(Q)$ for each non-singular matrix $Q$. It has the following
presentation:
$$
\bar{u} = u, $$
$$ u u^t = I_m = u^t u, $$
$$ u^t Q u Q^{-1} = I_m = Q u Q^{-1} u^t ,$$
where $u = (a_{ij})$.

As a matter of fact, it is more appropriate to use the notation
$A_o(Q)$ (and we will do so from now on)
for the compact matrix quantum group with the following sets of
relations (where $Q$ is positive):
$$ \bar{u} = u, $$
$$ u^t Q u Q^{-1} = I_m = Q u Q^{-1} u^t.$$
(Let $v = Q^{1/2} u Q^{-1/2}$. Then $v$ is a unitary matrix. Hence
the $C^*$-algebra $A$ exists. From this it is easy to see that $A_o(Q)$
is a compact matrix quantum group.)
This quantum group has all the properties
listed in \cite{W5} for the old $A_o(Q)$.
The old $A_o(Q)$ is the intersection of the
quantum groups $A_o(n)$ and the new $A_o(Q)$ defined
above. Moreover, if $Q$ is a real matrix, the new
$A_o(Q)$ is a compact quantum group of Kac type.

Finally, we note that the quantum group denoted by
$A_o(F)$ in \cite{Banica1} is the same as the quantum
group $B_u(Q)$ in \cite{W5',W7} with $Q = F^*$,
so it is different from the quantum group $A_o(Q)$ above
unless $F$ is the trivial matrix $I_n$.
\vvvv
{\bf Acknowlegdement.}
The author wishes to thank Alain Connes for several
helpful discussions and for his interest in this work.
He is also indebted to Marc Rieffel for his support,
which enabled the author finish writing up this paper. He
thanks T. Hodges, G. Nagy, A. Sheu, S.L. Woronowicz for their comments
during the AMS summer research conference on Quantization in July, 1996,
on which the author reported preliminary results of this paper.
The main results of this paper were obtained while the author was a visiting
member at the IHES during the year July, 1995-Aug, 1996.
He is grateful for the financial support of the IHES during this period.
He would like to thank the Director Professor
J.-P. Bourguignon and the staff of the IHES for their hospitality.
The author also wishes to thank the Department of
Mathematics at UC-Berkeley for its support and hospitality
while the author held an NSF Postdoctoral Fellowship there
during the final stage of this paper.

\vspace{1.5cm}
\hfill 
Sept, 1997

\end{document}